\documentclass[12pt]{amsart}

\usepackage{amssymb}

\usepackage{enumerate}

\usepackage{graphicx}

\makeatletter
\@namedef{subjclassname@2010}{%
  \textup{2010} Mathematics Subject Classification}
\makeatother

\newtheorem{theorem}{ Main Theorem}[section]
\newtheorem{thm}{Theorem}[section]

\newtheorem{conj}[thm]{Conjecture}
\newtheorem{pro}[thm]{Proof}

\theoremstyle{definition}

\newtheorem{rem}{Remark}
\newtheorem{exa}[thm]{Example}

\numberwithin{equation}{section}

\begin{document}


\baselineskip=17pt


\title[Is the quartic Diophantine equation $A^4+hB^4=C^4+hD^4$ solvable for any integer $h$? ]{Is the quartic Diophantine equation $A^4+hB^4=C^4+hD^4$ solvable for any integer $h$?}

\author[F. Izadi]{Farzali Izadi}
\address{Farzali Izadi \\
Department of Mathematics \\ Faculty of Science \\ Urmia University \\ Urmia 165-57153, Iran}
\email{f.izadi@urmia.ac.ir}

\author[M. Baghalaghdam]{Mehdi Baghalaghdam}
\address{Mehdi Baghalaghdam \\
Department of Mathematics\\ Faculty of Science \\ Azarbaijan Shahid Madani University\\Tabriz 53751-71379, Iran}
\email{mehdi.baghalaghdam@yahoo.com}

\date{}

\begin{abstract}
The Diophantine equation $A^4+hB^4=C^4+hD^4$, where $h$ is a fixed arbitrary positive integer, has been investigated by some authors. Currently, by computer search, the integer solutions of this equation are known for all positive integer values of $h \le 5000$ and $A, B, C, D \le 100000$, except for some numbers, while a solution of this Diophantine equation is not known for arbitrary positive integer values of $h$. Gerardin and Piezas found solutions
of this equation when $h$ is given by polynomials of degrees $5$ and $2$
respectively. Also Choudhry presented some new solutions of this
equation when $h$ is given by polynomials of degrees $2$, $3$, and $4$.

In this paper, by using the elliptic curves theory, we study this Diophantine equation, where $h$ is a fixed arbitrary rational number. We work out some solutions of the Diophantine equation for certain values of $h$, in particular for the values which has not already been found a solution in the range where $A, B, C, D \le 100000$ by computer search. Also we present some new parametric solutions for the Diophantine equation when $h$ is given by polynomials of degrees $3$, $4$. Finally We present two conjectures such that if one of them is correct, then we may solve the above Diophantine equation for arbitrary rational number $h$.
\end{abstract}



\maketitle

\section{Introduction}
\noindent The Diophantine equation (DE)
\begin{equation}\label{e8}
 A^4+hB^4=C^4+hD^4,
 \end{equation}
  where $h$ is a fixed arbitrary positive integer, has been investigated by some authors.
\\

\noindent The numerical solutions for
$75$ integer values of $h\leq 101$ was given by Choudhry \cite{A.C}. Then these solutions were
first extended by Piezas \cite{T.P} for all positive integer values of $h \leq 101$,
and finally by Tomita \cite{S.T} for all positive integer values of $h \le 1000$ except
$h = 967$. The lost solution for $h = 967$ was supplied by Bremner. Currently, by computer search, the small solutions of this DE are known for all positive integers $h \le 5000$, and $A, B, C, D \le 100000$ except for
\\

\noindent $h= 1198, 1787, 1987$,\\

\noindent $2459, 2572, 2711, 2797, 2971$,\\

\noindent $3086, 3193, 3307, 3319, 3334, 3347, 3571, 3622, 3623, 3628, 3644$,\\
 $3646, 3742, 3814, 3818, 3851, 3868, 3907, 3943, 3980$,\\

\noindent  $4003, 4006, 4007, 4051, 4054, 4099, 4231, 4252, 4358, 4406, 4414$,\\
 $4418, 4478, 4519, 4574, 4583, 4630, 4643, 4684, 4870, 4955, 4999.$
\\

\noindent We will work out some of these cases.
\\

\noindent Gerardin and Piezas found solutions
of this equation when $h$ is given by polynomials of degrees $5$, and $2$,
respectively, see  \cite{T.P} and \cite{S.T}. Also Choudhry presented several new solutions of this
equation when $h$ is given by polynomials of degrees $2$, $3$, and $4$, see \cite{A.C2}.
\\

\noindent In this paper, we used  elliptic curves theory to study the DE \eqref{e8}. \\

\section{ The first method for solving the DE $A^4+hB^4=C^4+hD^4$}

\noindent Our main result in this section is the following:

\begin{theorem} Consider the DE \eqref{e8},
where $h$ is a fixed arbitrary rational number.\\ Then there exists a cubic elliptic curve of the form

\noindent $E(h):Y^2=X^3+FX^2+GX+H$, where the coefficients $F$, $G$, and $H$, are all functions of $h$.
If the elliptic curve $E(h)$ or its counterpart $E(h)_t$ resulting from $E(h)$ by switching $h$ to $ht^4$
has positive rank, depending on the value of $h$ and an appropriate rational number $t$, then the DE \eqref{e8} has infinitely many integers solutions.
By taking $h=\frac{v}{u}$, this also solves DE of the form $uA^4+vB^4=uC^4+vD^4$ for appropriate integer values of $u$ and $v$.
\end{theorem}

\begin{pro} Let: $A=m-q$, $B=m+p$, $C=m+q$, and $D=m-p$, where all variables are rational numbers. By substituting these variables in the DE \eqref{e8} we get
\begin{equation}
-8m^3q-8mq^3+8hm^3p+8hmp^3=0.
\end{equation}

\noindent Then after some simplifications and clearing the case $m=0$ we obtain
\begin{equation}\label{e7}
m^2(hp-q)=-hp^3+q^3.
\end{equation}

\noindent We may assume that $hp-q=1$ and $m^2=-hp^3+q^3$.

\noindent By plugging $q=hp-1$
 into the equation \eqref{e7} and some simplifications we obtain the equation

\begin{equation}\label{e1}
m^2=(h^3-h)p^3-(3h^2)p^2+(3h)p-1.
\end{equation}
\\
\noindent By multiplying both sides of this equation in
 $(h^3-h)^2$ and letting

\begin{equation} \label{e2}
  X=(h^3-h)p \hspace{1cm} Y=(h^3-h)m,
\end{equation}
\\
\noindent we get  the elliptic curve

\begin{equation}\label{e4}
 Y^2=X^3-(3h^2)X^2+(3h(h^3-h))X-(h^3-h)^2.
 \end{equation}
\\
\noindent By Letting $X=Z+h^2$ in \eqref{e4}, we get the simple elliptic curve

\begin{equation}\label{e78}
E(h): Y^2=Z^3-(3h^2)Z-(h^4+h^2).
\end{equation}

\noindent If for  a given $h$, the above elliptic curve $E(h)$ or its counterpart $E(h)_t$ resulting from $E(h)$ by switching $h$ to $ht^4$ has positive rank, then by calculating $m$, $p$, $q$, $A$, $B$, $C$, $D$, from the relations \eqref{e2}, $q=hp-1$, $A=m-q$, $B=m+p$, $C=m+q$, $D=m-p$, after some simplifications  and canceling  the  denominators of $A$, $B$, $C$, $D$, we obtain infinitely many integer solutions for the DE \eqref{e8}. Now the proof of the main theorem is completed.
\end{pro}

\noindent Although, we were able to find an appropriate $t$ such that $E(h)_t$ has positive rank in the case of rank zero $E(h)$ for many values of $h$,
the proof for arbitrary $h$ seems to be difficult at this point. For this reason, we state it as a conjecture.
\\

\begin{conj}
Let $h$ be an arbitrary fixed rational number. Then there  exists  at least a rational number $t$ such that the rank of the elliptic curve
\begin{equation}
E(h)_t: Y^2=Z^3-(3h^2t^8)Z-(h^4t^{16}+h^2t^8),
 \end{equation}

\noindent is positive.
\end{conj}

\begin{rem}
 If $h$ is a large number, for example $h=7000$, we may divide $h$ by $10^4$  to get an elliptic curve of positive rank  with small  coefficients, then solve the DE \eqref{e8} for $h'=\frac{7}{10}$ and finally get a solution for the main case of $h=7000$,
by multiplying  both sides of the DE \eqref{e8} by an appropriate number. As another example, if $h=9317=7.11^3$,
we may first solve the DE \eqref{e8} for $h'=\frac{7}{11}$.
\end{rem}

\begin{rem} Note that by substituting the relation $p=\frac{q+1}{h}$, in the equation
 $m^2=-hp^3+q^3$, multiplying both side of the equation by $(\frac{h^2-1}{h^2})^2$, and letting

\begin{equation}\label{e18}
 X'=(\frac{h^2-1}{h^2})q  \hspace{1cm} Y'=(\frac{h^2-1}{h^2})m,
\end{equation}
\\
\noindent we get another elliptic curve

\begin{equation}\label{e21}
Y'^2=X'^3-\frac{3}{h^2}X'^2-\frac{3(h^2-1)}{h^4}X'-(\frac{h^2-1}{h^3})^2.
\end{equation}
\\

\noindent Now, if we set $Y'=\frac{Y}{h^3}$, $X'=\frac{X+1-h^2}{h^2}$, then the elliptic curve \eqref{e21}, transforms to the elliptic curve \eqref{e4}. This means that two elliptic curves \eqref{e21} and \eqref{e4} are isomorphic.
\end{rem}

\section{ Application to examples}

\subsection{Example: $A^4+B^4=C^4+D^4$}\noindent

\noindent i.e., sums of two biquadrates in two different ways.
\\

\noindent $h=16$, here $h=1$ replaced by $h=2^4$.\\

\noindent $E(16)$: $Y^2=X^3-768X^2+195840X-16646400$.\\

\noindent Rank=1.\\

\noindent Generator: $P=(X,Y)=(340,680)$.\\

\noindent Points:  $2P=(313,-275)$, $3P=(\frac{995860}{729},\frac{-727724440}{19683})$,
 $4P=(\frac{123577441}{302500},\frac{305200800239}{166375000})$.\\

\noindent $(p',m',q')=(\frac{313}{4080},\frac{-55}{816},\frac{58}{255})$, \\
 $(p'',m'',q'')=(\frac{2929}{8748},\frac{-1070183}{118098},\frac{9529}{2187})$,\\
$(p''',m''',q''')=(\frac{123577441}{1234200000},\frac{305200800239}{678810000000},\frac{46439941}{77137500})$.\\

\noindent Solutions:\\

$1203^4+76^4=653^4+1176^4$,\\

$1584749^4+2061283^4=555617^4+2219449^4$,\\

$103470680561^4+746336785578^4=713872281039^4+474466415378^4$.\\

\subsection{Example: $A^4+206B^4=C^4+206D^4$}\noindent
\\

\noindent $h=\frac{103}{8}$.\\

\noindent $E(\frac{103}{8})$: $Y^2=X^3-\frac{31827}{64}X^2+\frac{335615715}{4096}X-\frac{1179689238225}{262144}$.\\

\noindent Rank=1.\\

\noindent Generator: $P=(X,Y)=(\frac{2131205}{32},\frac{8767168835}{512})$.\\

\noindent $(p,m,q)=(\frac{6819856}{217227},\frac{1753433767}{217227},\frac{850373}{2109})$.\\

\noindent Solution:\\

$3331690696^4+206.(1760253623)^4=3682044372^4+206.(1746613911)^4$.\\

\begin{rem}By searching, Noam Elkies found the smallest solution to this DE as $A, B, C, D=3923, 1084, 4747, 506.$
\end{rem}

\begin{rem} No computer research has come up with a solution for the following equations in the range of
$A, B, C, D \le 100000$, see \cite{S.T}.
\end{rem}

\subsection{Example: $A^4+2572B^4=C^4+2572D^4$}

\noindent $h=2572$.\\

\noindent $E(2572)$: $Y^2=X^3-(3.2572^2)X^2+(3.2572.(2572^3-2572))X-(2572^3-2572)^2$.\\

\noindent Rank=2.\\
\noindent Generators: $P_1=(X,Y)=(\frac{60035809}{9},\frac{302757191}{27})$,\\

 \noindent $P_2=(X',Y')=(\frac{3435573760731933430513}{381659437643236},\frac{27488556048550361767336062809879}{7456139229698648679016})$. \\

\noindent $(p,m,q)=(\frac{23333}{59513508},\frac{117667}{178540524},\frac{194}{23139})$,
 \\

\noindent $(p',m',q')=(\frac{1558040235953533}{2944884220855208976},\frac{12466120460409195830562539}{57531570296354773207287456},\frac{413061923023825}{1144978312929708})$.
\\

\noindent Solutions:\\

$1379237^4+2572.(187666)^4=1614571^4+2572.(47668)^4$,
\\

$8288946070402543055294861^4+2572.(12496558499611049062325037)^4=\\
33221186991220934716419939^4+2572.(12435682421207342598800041)^4$.
\\

\subsection{Example: $A^4+967B^4=C^4+967D^4$}\noindent

\noindent $h=967$.\\

\noindent $E(967)$: $Y^2=X^3-(3.967^2)X^2+(3.967.(967^3-967))X-(967^3-967)^2$.\\

\noindent Rank=1.\\

\noindent Generator: $P=(X,Y)=(\frac{238501273696}{245025},\frac{900632541139856}{121287375} )$.\\

\noindent $(p,m,q)=(\frac{2129475658}{1978205172075},\frac{8041361974463}{979211560177125},\frac{83761933}{2045713725})$.\\

\noindent Solution:\\

$32052543684982^4+967.(9095452425173)^4=\\
48135267633908^4+967.(6987271523753)^4$.
\\

\subsection{Example: $A^4+2797B^4=C^4+2797D^4$}

\noindent $h=2797$.\\

\noindent $E(2797)$: $Y^2=X^3-(3.2797^2)X^2+(3.2797.(2797^3-2797))X-(2797^3-2797)^2$.\\

\noindent Rank=1.\\

\noindent Generator: $P=(X,Y)=(\frac{18256234369}{2304},\frac{3411597220289}{110592})$.\\

\noindent $(p,m,q)=(\frac{6527077}{18024671232},\frac{1219734437}{865184219136},\frac{231563137}{18024671232})$.\\

\noindent Solution:\\

$9895296139^4+2797.(1533034133)^4=12334765013^4+2797.(906434741)^4$.\\

\subsection{Example: $A^4+4999B^4=C^4+4999D^4$}

\noindent $h=4999$.\\

\noindent $E(4999)$: $Y^2=X^3-(3.4999^2)X^2+(3.4999.(4999^3-4999))X-(4999^3-4999)^2$.\\

\noindent Rank=1.\\

\noindent Generator:$P=(X,Y)=(\frac{38932053386017900293094583125}{1502165941669975655844},\frac{51963991529347119364735376770810745620625}{58220625445784716642962064124328})$.\\

\noindent  $p=\frac{62291285417628640468951333}{300252952455649800809709213504}$,\\

\noindent $m=\frac{83142386446955390983576602833297192993}{11637139545633456295150723664563706629248},$\\

\noindent $q=\frac{2228682405896333845684837}{60062603011732306623266496}$.\\

\noindent Solution:\\

\noindent $A=348665208625932834629908938859838853613$,\\
$B=85556658729553179445421813716725247139$,\\
$C=514949981519843616597062144526433239599$,\\
$D=80728114164357602521731391949869138847$.\\
\\

\subsection{Example: $A^4+2459B^4=C^4+2459D^4$}

\noindent $h=2459$.\\

\noindent $E(2459)$: $Y^2=X^3-(3.2459^2)X^2+(3.2459.(2459^3-2459))X-(2459^3-2459)^2$.\\

\noindent Rank=1.\\

\noindent  Generator: $P=(X,Y)$, where
\\

\noindent $X=\frac{2455940168334175449299068876662469864}{403764781843031846693075441721}$,\\

\noindent $Y=\frac{775339319798703416232888693955985044070341765700696}{256562189232730518019448407676170227655852269}$.\\

\noindent  $p=\frac{249790497186144777186642481352977}{610607423109489416276237594383042485}$,\\

\noindent $m=\frac{78858759133309948762498850076890260788277234103}{387995150343819871443041535202557092730530548624665},$\\

\noindent $q=\frac{1475156352680191470401084694562}{248315340833464585716241396658415}$.
 \\

\noindent Solution:\\

\noindent $A=2226087479458719030508635008690035436036778215959$,\\
$B=237581856564140327136761830581698727005176625756$,\\
$C=2383804997725338928033632708843815957613332684165$,\\
$D=79864338297520429611764130427918205428622157550$.\\

\section{The second method for solving the DE $A^4+hB^4=C^4+hD^4$}

\noindent In this section, we wish to look at the equation from a different perspective. To this end we take $X=Z^2+h^2$  \eqref{e4} to get the six
degree curve

\begin{equation}
Y^2=Z^6-3h^2Z^2-(h^4+h^2).
\end{equation}

\noindent This curve can be considered as an quatic elliptic curve of $(h,Y)$ letting $Z$ as a parameter, i.e.,

\begin{equation}\label{e100}
Y^2=-h^4-(3Z^2+1)h^2+Z^6.
\end{equation}

\noindent Next we use theorem $2.17.$ of  \cite{L.W} to transform this quartic to a cubic elliptic curve of the form

\begin{equation}\label{e500}
E'(Z): Y'^2=X'^3-(3Z^2+1)X'^2+(4Z^6)X'-(12Z^8+4Z^6).
\end{equation}

\noindent with the inverse transformation

\begin{equation}\label{e600}
h=\frac{2Z^3(X'-(3Z^2+1))}{Y'} \hspace{1cm} Y=-Z^3+\frac{h^2X'}{2Z^3}.
\end{equation}

\noindent Now by taking an appropriate rational value for $Z$ such that the rank of the elliptic curve \eqref{e500} is positive, we obtain an infinitely many rational points on \eqref{e500} and consequently an infinite set of rational values for $h$ (also for $Y$)), which is denoted by $H(Z)$ . Then for every $h$ obtained in this way, we can find a solution for the \eqref{e4} and finally a solution for the main DE by using all the necessary transformations namely

\noindent $X=Z^2+h^2$, $m=\frac{Y}{h^3-h}$,  $p=\frac{X}{h^3-h}$, $h=\frac{q+1}{p}$, $A=m-q$, $B=m+p$, $C=m+q$, $D=m-p$.
\\

\noindent To get infinitely many solutions one can use the Richmond method \cite{H.R}. The following examples clarify this idea better

\noindent $Z_1=3$.

\noindent $E'(3)$: $Y'^2=X'^3-28X'^2+2916X'-81648$.

\noindent Rank=1.

\noindent Generator: $P=(X',Y')=(108,1080)$.

\noindent  $H(Z_1)=\{4,\frac{2^3.3^3}{197},\frac{-2^2.251.395449}{11.13.61.653}, \cdots \}$.\\

\noindent $Z_2=4$.

\noindent $E'(4)$: $Y'^2=X'^3-49X'^2+16384X'-802816$.

 \noindent Rank=1.

\noindent Generator: $P=(X',Y')=(202,2958)$.

\noindent $H(Z_2)=\{\frac{2^6.3}{29},\frac{2^9.3.29}{61121},\frac{-2^6.3^2.19.3571.18131}{5.29.97.7746413},\cdots \}$.\\

\noindent $Z_3=6$.

\noindent $E'(6)$: $Y'^2=X'^3-109X'^2+186624X'-20342016$.

\noindent Rank=1.

\noindent  Generator: $P=(X',Y')=(\frac{621}{4},\frac{-24975}{8})$.

\noindent $H(Z_3)=\{\frac{-2^2.5.317}{3^3.37},\frac{2^3.5^4.408841}{17.37^2.12757}$,$\frac{-2^2.5.10193.249587558933}{3^3.7.37.5101.181680953},\cdots \}$.\\

\noindent Having seen the examples, the natural question arises:

\noindent  Does the set of natural numbers $\mathbb{N}$ contained in $\bigcup_ { t \in \mathbb{Q^*}}t^4 (\bigcup_{Z \in \Omega}H(Z))$?,

\noindent where

\noindent $\Omega=\{Z \in \mathbb{Q} \mid E'(Z) \hspace{.1cm}has \hspace{.1cm} positive \hspace{.1cm}rank \}.$

\noindent We state this as the second conjecture:

\begin{conj}
With the above notations one has $\mathbb{N} \subset \bigcup_ { t \in \mathbb{Q^*}}t^4 (\bigcup_{Z \in \Omega}H(Z)) . $
\end{conj}

\section{Application to examples}\noindent
\\

\noindent Now we are going to work out some examples.
\subsection {Example: $h=108$} \noindent

\noindent $Z=\frac{5}{3}$.

\noindent $E'(\frac{5}{3})$: $Y'^2=X'^3-\frac{-28}{3}X'^2+\frac{62500}{729}X'-\frac{1750000}{2187}$.

\noindent Rank=1.

\noindent Generator: $P=(X',Y')=(\frac{2500}{81},\frac{109000}{729})$.

\noindent $(m,p,q)=(\frac{5}{4},\frac{123}{28},\frac{34}{7})$.

\noindent $h=\frac{4}{3}$.

\noindent Solution:

$303^4+108(158)^4=513^4+108(88)^4$.\\

\noindent Note: $108=4.27$.

\subsection {Example: $h=492$} \noindent

\noindent $Z=\frac{4}{3}.$

\noindent $E'(\frac{4}{3})$: $Y'^2=X'^3-\frac{19}{3}X'^2+\frac{16384}{729}X'-\frac{311296}{2187}$.

\noindent Rank=1.

\noindent Generator: $P=(X',Y')=(\frac{586}{81},\frac{5986}{729})$.

\noindent $(m,p,q)=(\frac{56908}{11033},\frac{-238251}{44132},\frac{-42025}{11033})$.

\noindent $h=\frac{64}{123}$.

\noindent Solution:

$42476^4+492(395732)^4=1863532^4+492(59532)^4$.\\

\noindent Note: $492=123.4$.

\subsection {Example $h=12256974$} \noindent
\\
\noindent $Z=\frac{3}{2}$.

\noindent $E'(\frac{3}{2})$: $Y'^2=X'^3-\frac{31}{4}X'^2+\frac{729}{16}X'-\frac{22599}{64}$.

\noindent Rank=1.

\noindent Generator: $P=(X',Y')=(\frac{135}{4},\frac{351}{2})$.

\noindent $2P=(\frac{665}{64},\frac{10309}{512})$.

\noindent $h=\frac{54}{61}$.

\noindent $(m,p,q)=(\frac{145851}{12880},\frac{-306037}{19320},\frac{-48373}{3220})$.

\noindent Solutions:

$62099769^4+12256974(174521)^4=8718303^4+12256974(1049627)^4$.\\

\noindent Note: $54.61^3=12256974$.

\noindent For the point 3P, we obtain $h=\frac{-805}{3977}$, then as two more examples, we obtain solutions for the  two cases $h'=3977.805^3$, $h''=3977^3.805$.

\section{Parametric solutions of $A^4+hB^4=C^4+hD^4$}
\noindent Let: $A=m-q$, $B=m+p$, $C=m+q$, and, $D=m-p$, where all variables are rational numbers. By substituting these variables in the DE \eqref{e8}, and some simplification, we get

\begin{equation}\label{a}
p(hm^2+hp^2)=q(m^2+q^2).
\end{equation}

\noindent We may assume that

 \begin{equation}\label{b}
p=m^2+q^2,
\end{equation}
and
\begin{equation}\label{c}
q=h(m^2+p^2).
\end{equation}

 \noindent By substituting $p=m^2+q^2$, in the relation $h=\frac{q}{m^2+p^2}$, we get:

\begin{equation}\label{b}
 h=\frac{q}{m^2+(m^2+q^2)^2}=\frac{q}{m^2+m^4+q^4+2m^2q^2}.
 \end{equation}
\\

\noindent Then using the reciprocal of $h$, we conclude that
\begin{equation}
 h=\frac{m^2+(m^2+q^2)^2}{q}=\frac{m^2+m^4+q^4+2m^2q^2}{q},
\end{equation}
\noindent  the DE \eqref{e8} has a parametric solution:

\noindent $A=m+m^2+q^2$,\\
$B=m-q$,\\
$C=m-m^2-q^2$,\\
$D=m+q$.\\

\begin{exa}
$m=kq$;\\
$h=k^2q+(k^2+1)^2q^3;$\\
$A=k+k^2q+q$;\\
$B=k-1$;\\
$C=k-k^2q-q$;\\
$D=k+1$;\\

\begin{rem}
The above example provides parametric solutions for the DE \eqref{e8} when $h$ is given by polynomials of the degrees $3$, and $4$ as follows.
\end{rem}

\noindent $m=q$;\\
$h=q+4q^3;$\\
$A=1+2q$;\\
$B=0$;\\
$C=1-2q$;\\
$D=2$;\\

\noindent If we let $2q=p$, this recovers  the third and the second parametric solutions of the table $1$ in \cite{A.C2}. Of course by replacing $p$ with $2p$ in  the second parametric solution of the table $1$ in \cite{A.C2}, and dividing both sides of this equation  by $16$, we get the
third parametric solution of the table $1$.!\\
\\
$h=8p(p^2+1);$\\
$A=p+1$;\\
$B=0$;\\
$C=p-1$;\\
$D=1$;\\

\noindent $m=2q$;\\
$h=4q+25q^3;$\\
$A=2+5q$;\\
$B=1$;\\
$C=2-5q$;\\
$D=3$;\\

\noindent $m=3q$;\\
$h=9q+100q^3;$\\
$A=3+10q$;\\
$B=2$;\\
$C=3-10q$;\\
$D=4$;\\

\noindent $q=2$, $m=2k$;\\
$h=8k^4+18k^2+8$;\\
$A=2k^2+k+2$;\\
$B=k-1$;\\
$C=2k^2-k+2$;\\
$D=k+1$;\\

\noindent $q=3$, $m=3k$;\\
$h=27k^4+57k^2+27$;\\
$A=3k^2+k+3$;\\
$B=k-1$;\\
$C=3k^2-k+3$;\\
$D=k+1$;\\
\end{exa}

\begin{exa}\noindent

 \noindent $h=8k^3p+512k^7p^3+512k^3p^3+1024k^5p^3$;\\
$A=8k^3p+8kp-k$;\\
$B=k+1$;\\
$C=8k^3p+8kp+k$;\\
$D=k-1$;\\

\noindent $k=\frac{1}{2}$;\\
$h=100p^3+p$;\\
$A=10p-1$;\\
$B=3$;\\
$C=10p+1$;\\
$D=1$;\\
\end{exa}

\begin{exa}\noindent \\
$h=n(p^4+(n^2+2)p^2+1)$;\\
$A=p^2+np+1$;\\
$B=p+1$;\\
$C=p^2-np+1$;\\
$D=p-1$;\\
\end{exa}
\begin{rem}
\noindent The case $n=1$, recovers the last parametric solution of the table 1 in \cite{A.C2}.
\end{rem}

\begin{exa}
\noindent $q=1$;\\
$h=m^4+3m^2+1;$\\
$A=m+m^2+1$;\\
$B=m-1$;\\
$C=m-m^2-1$;\\
$D=m+1$;\\
\end{exa}
\noindent Again this recovers the last parametric solution of the table 1 in \cite{A.C2}.

\begin{exa}\noindent \\
$h=(p^2+2)(p^2+4)$;\\
$A=p^2+p+2$;\\
$B=p+1$;\\
$C=p^2-p+2$;\\
$D=p-1$;\\
\end{exa}
\begin{exa}\noindent

\noindent$h=512m^4+1032m^2+512$;\\
$A=8m^2-m+8;$\\
$B=m+1$;\\
$C=8m^2+m+8$;\\
$D=m-1$;\\
\end{exa}

\begin{exa}
\noindent $q=m^2$;\\
$h=1+m^2+m^6+2m^4;$\\
$A=1+m+m^3$;\\
$B=1-m$;\\
$C=1-m-m^3$;\\
$D=1+m$;\\
\end{exa}

\noindent\begin{exa}
By taking $Z=h^2+1$, in the elliptic curve \eqref{e78} (or $X=2h^2+1$ in the \eqref{e4}), we get
\begin{equation}
Y^2=h^6-h^4-h^2+1=(h^2-1)(h^4-1)=(h^2-1)^2(h^2+1).
\end{equation}

\noindent letting $h^2+1=t^2$, yields $h=\frac{r^2-1}{2r}$ and $t=\frac{-(r^2+1)}{2r}$, then we get
\\

\noindent $Y=(h^2-1)(-t)=\frac{(r^6-5r^4-5r^2+1)}{8r^3}$, and $X=\frac{r^4+1}{2r^2}$.\\

\noindent Finally by calculating $m, p, q, A, B, C, D$, from the above relations, we obtain the parametric solution as follows:\\
$h=8r^3(r^2-1)$,\\
$A=4r(5r^4-1)$,\\
$B=(r^4+6r^3+6r^2+6r+1)(r-1)^2$,\\
$C=4r^3(r^4-5)$,\\
$D=(r^4-6r^3+6r^2-6r+1)(r+1)^2$.\\

\noindent Note that $h$ is given by polynomial of degree $5$.

\end{exa}

\noindent The Sage software and Denis Simon's ellrank code were used for calculating  the rank of the elliptic curves, (see \cite{S.A}).
\\

\begin{center}\textbf{Acknowledgements}

\noindent We are very grateful to Professor Allen MacLeod for the rank and generator computations of the elliptic curves with big generators. Finally the second author would like to present this work to his parents and his wife.
\end{center}

\end{document}